\documentclass{amsart}

\usepackage{amssymb,mathtools}
\usepackage{tikz}
\usepackage[colorlinks]{hyperref}

\newtheorem{theorem}{Theorem}[section]
\newtheorem{proposition}[theorem]{Proposition}
\newtheorem{lemma}[theorem]{Lemma}
\newtheorem{corollary}[theorem]{Corollary}
\theoremstyle{definition}

\newtheorem{example}[theorem]{Example}
\theoremstyle{remark}
\newtheorem{remark}[theorem]{Remark}

\DeclareMathOperator{\pdim}{pdim}
\DeclareMathOperator{\idim}{idim}
\DeclareMathOperator{\fdim}{fdim}
\DeclareMathOperator{\Gpdim}{Gpdim}
\DeclareMathOperator{\Ext}{Ext}
\DeclareMathOperator{\Hom}{Hom}
\DeclareMathOperator{\domdim}{domdim}

\DeclareMathOperator{\grade}{grade}
\DeclareMathOperator{\depthr}{depth}

\DeclareMathOperator{\Spec}{Spec}
\DeclareMathOperator{\Max}{Max}
\DeclareMathOperator{\soc}{soc}
\DeclareMathOperator{\Sim}{Sim}
\DeclareMathOperator{\ann}{ann}
\DeclareMathOperator{\Tor}{Tor}
\DeclareMathOperator{\op}{op}
\renewcommand{\mod}{\operatorname{mod}}
\newcommand{\lto}{\longrightarrow}

\title[An Auslander--Bridger formula]
{An Auslander--Bridger formula for Auslander--Gorenstein algebras}
\author{Dawei Shen}
\address{School of Mathematics and Statistics, 
    Henan University, 
    Kaifeng, 
    Henan 475004, 
    P. R. China}
\email{sdw12345@mail.ustc.edu.cn}
\date{\today}

\begin{document}
\begin{abstract}
For an Auslander--Gorenstein Noether algebra 
finite over a commutative Noetherian ring, we express the Gorenstein
projective dimension of a nonzero finitely generated module as a supremum
of the differences of depths measured with respect to simple modules. Among
Iwanaga--Gorenstein Artin algebras, validity of the formula characterizes
the Auslander condition. For Iwanaga--Gorenstein Artin algebras of positive
self-injective dimension, the number of simple modules of maximal Gorenstein
projective dimension is bounded
by the number of nonprojective indecomposable injectives, with equality
precisely for minimal Auslander--Gorenstein algebras.
\end{abstract}
\maketitle

\section{Introduction}
The Auslander--Buchsbaum formula expresses finite projective dimension
over a commutative local Noetherian ring as a difference of depths \cite{AusBuch}. 
Auslander and Bridger extend this identity
to finite G-dimension \cite{AB}. In both formulas, depth is measured
by the first nonvanishing Ext group from the residue field.
In noncommutative versions, an analogous role is played by simple
modules or by dominant dimension. 
Wu and Zhang obtain an Auslander--Buchsbaum formula
in an Artin--Schelter Gorenstein setting \cite{WZ}.
For connected minimal Auslander--Gorenstein Artin algebras,
Chen, Iyengar and Marczinzik relate finite projective dimension
to dominant dimension for modules that are not projective-injective
\cite{ChIM}. Cruz and Marczinzik give another proof of this formula for
higher Auslander algebras \cite{CM}.

Auslander--Gorenstein rings provide a noncommutative analog of
commutative Gorenstein rings. Over Noether algebras, the
classification of indecomposable injectives involves simple modules
over central localizations at all prime ideals \cite{KN}. In particular,
terms associated with nonmaximal prime ideals need not contain simple
modules over the original algebra. An extension of the classical
depth formulas therefore requires an explanation of how much of this
injective structure can be detected by simple modules 
over the original algebra.

Let $R$ be a commutative Noetherian ring and let $A$ be a unital
$R$-algebra that is finitely generated as an $R$-module.
Such an algebra is called a Noether $R$-algebra. 
We write $\mod A$ for the category of finitely generated $A$-modules.

An algebra $A$ is called Iwanaga--Gorenstein if its self-injective
dimensions are finite on both sides.
Let $I^\bullet(A)$ be a minimal injective coresolution
of the regular module. The Auslander condition requires that
the flat dimension of $I^q(A)$ be at most $q$ for every $q\ge0$.
An Iwanaga--Gorenstein algebra satisfying this condition is
said to be Auslander--Gorenstein.

When $A$ is Auslander--Gorenstein, the main result shows that simple
modules suffice to compute the Gorenstein projective dimension of every
nonzero finitely generated $A$-module $M$.
For each simple module $S$, 
the formula compares the $S$-depth of $A$ with the $S$-depth of $M$.
The supremum of these differences recovers the Gorenstein projective
dimension of $M$.

Following \cite{WZ}, 
for a simple $A$-module $S$ and a finitely generated $A$-module $M$, 
we set 
\[
 \depthr_S M=\inf\{i\ge0\mid\Ext_A^i(S,M)\ne0\}.
\]
We adopt the convention that the infimum of the empty set is $\infty$.

The following Auslander--Bridger formula is our main result.

\begin{theorem}\label{intro:formula}
Let $A$ be an Auslander--Gorenstein Noether algebra.
For any nonzero $M$ in $\mod A$, we have
\[
 \Gpdim_A M=\sup_{S\in\Sim A}
 \{\depthr_S A-\depthr_S M\}.
\]
In particular, $M$ is Gorenstein projective if and only if
\[\depthr_S M\ge\depthr_S A\]
for every simple module $S$.
\end{theorem}

Over a commutative Gorenstein ring of finite Krull dimension,
the last assertion recovers the characterization of finitely generated
Gorenstein projective modules as maximal Cohen--Macaulay modules.

The key step is to detect an injective summand in the minimal
injective coresolution of $M$.
Tor nonvanishing and vanishing between injectives associated with
distinct central primes ensure that this summand is the injective
envelope of a simple $A$-module.
This detection requires only the Iwanaga--Gorenstein hypothesis.
The Auslander condition converts the detection bound into
the Auslander--Bridger formula.

In the Artin setting, the formula also follows from the known
characterization of Gorenstein projective dimension by flat-dimension
bounds on the terms of minimal injective coresolutions.
Combined with the simple-module characterization of the
Auslander condition, it gives the following converse.

\begin{corollary}\label{intro:artin}
Let $A$ be an Iwanaga--Gorenstein Artin algebra.
The following are equivalent.
\begin{enumerate}
\item $A$ is Auslander--Gorenstein;
\item The Auslander--Bridger formula holds for every nonzero $M\in\mod A$.
\end{enumerate}
\end{corollary}

We next apply the criterion of Li, Marczinzik and Zhang \cite{LMZ}
to obtain a characterization of minimal Auslander--Gorenstein Artin
algebras.

For an Artin algebra $A$, the defect $\operatorname{def} A$ counts the
isomorphism classes of nonprojective indecomposable injective modules.
Recall that an Artin algebra $A$ is minimal Auslander--Gorenstein
if its self-injective dimension is finite and bounded by its dominant
dimension.

\begin{theorem}\label{intro:defect}
Let $A$ be an Iwanaga--Gorenstein Artin algebra with positive
self-injective dimension $d$.
Then
\[
 \#\{S\in\Sim A\mid\Gpdim_A S=d\}\le\operatorname{def} A.
\]
Equality holds if and only if $A$ is minimal Auslander--Gorenstein.
\end{theorem}

The paper is organized as follows. Section~2 recalls the classification
of indecomposable injectives over Noether algebras and proves the
Tor vanishing needed for the detection argument.
In Section~3, we establish the main formula and its converse for Artin
algebras. Section~4 relates the number of simple modules of maximal Gorenstein
projective dimension to the defect of an Artin algebra. Section~5 compares the
formula with known identities.

Throughout this paper, all modules are left modules unless otherwise stated.

\section{Preliminaries on Noether algebras}
Let $A$ be a Noether algebra.
Write $\Sim A$ for the set of isomorphism classes of
simple $A$-modules.
For any $A$-module $X$, write $I^\bullet(X)$ for its minimal
injective coresolution, and write
$I_A(X)$ for its injective envelope.

For a simple module $S$ and $M\in\mod A$, set
\[
 \depthr_S M=\inf\{i\ge0\mid \Ext_A^i(S,M)\ne0\}.
\]
We write $\grade_A S=\depthr_S A$. In the formulas below,
a finite number minus $\infty$ is interpreted as $-\infty$.

\begin{lemma}\label{lem:occurrence}
Let $S$ be a simple $A$-module and let $M\in\mod A$.
For every $i\ge0$, $\Ext_A^i(S,M)$ is nonzero if and only if
$I_A(S)$ is a direct summand of $I^i(M)$.
\end{lemma}
\begin{proof}
The kernel of each differential of $I^\bullet(M)$ is essential in its source,
so every map from $S$ into $I^i(M)$ has image in that kernel.
Then
\[\Ext_A^i(S,M)\cong\Hom_A(S,I^i(M)).\]
A nonzero map from $S$ to $I^i(M)$
extends to a map from $I_A(S)$ by the injectivity of $I^i(M)$.
Its kernel meets the essential submodule $S$ trivially,
so the extension is injective. It splits because $I_A(S)$ is injective.
The converse follows by restricting a summand embedding to $S$.
\end{proof}

For a Noether algebra $A$, denote by
$\operatorname{ind-Inj} A$ the set of isomorphism classes of
indecomposable injective $A$-modules. 

The following classification is established in \cite[Section~2]{KN}.

\begin{lemma}\label{lem:injective-classification}
Let $A$ be a Noether $R$-algebra. There is a bijection
\[
 \begin{aligned}
 \bigsqcup_{\mathfrak p\in\Spec R}\Sim A_{\mathfrak p}
 &\xrightarrow{\ \sim\ }
 \operatorname{ind-Inj} A,\\
 (\mathfrak p,S)&\longmapsto I_{A_{\mathfrak p}}(S).
 \end{aligned}
\]
Furthermore, each $r\in R\setminus\mathfrak p$
acts invertibly on $I_{A_{\mathfrak p}}(S)$,
and each element of $I_{A_{\mathfrak p}}(S)$ is annihilated by a power
of $\mathfrak p$.
\end{lemma}

We need the following.

\begin{lemma}\label{lem:closed-fiber}
Let $A$ be a Noether $R$-algebra and let $\mathfrak p,\mathfrak q$ be
distinct prime ideals of $R$.
For a simple $A_{\mathfrak p}$-module $S$ and a simple
$A_{\mathfrak q}^{\mathrm{op}}$-module $T$, one has
\[
 \Tor_i^A(I_{A_{\mathfrak q}^{\mathrm{op}}}(T),I_{A_{\mathfrak p}}(S))=0
\]
for every $i\ge0$.
\end{lemma}
\begin{proof}
Let $X$ be an
$A$-module and $Y$ an $A^{\mathrm{op}}$-module. For $r\in R$, the
endomorphisms of $\Tor_i^A(Y,X)$ induced by
multiplication by $r$ on $X$ and on $Y$ agree. Indeed, on a
projective resolution $P_\bullet$ of $X$, multiplication by 
$r$ is a chain lift of multiplication on $X$,
and the balancing relation on $Y\otimes_A P_\bullet$ identifies
the two actions.

If $r$ acts invertibly on one argument, it acts invertibly on
Tor by functoriality. If it acts locally nilpotently on $Y$,
write
\[
 Y=\operatorname*{colim}_{e\ge1}Y_e,\,Y_e=\{y\in Y:r^e y=0\}.
\]
Centrality makes each $Y_e$ an $A^{\mathrm{op}}$-submodule. Since Tor
commutes with filtered colimits and $r^e$ kills
$\Tor_i^A(Y_e,X)$, the induced action on
$\Tor_i^A(Y,X)$ is locally nilpotent. The same
argument applies to the other variable. An invertible locally
nilpotent endomorphism can occur only on the zero module.

If $\mathfrak p\setminus\mathfrak q$ is nonempty, choose
$r$ in this difference. By Lemma~\ref{lem:injective-classification},
it acts locally nilpotently on
$I_{A_{\mathfrak p}}(S)$ and invertibly on $I_{A_{\mathfrak q}^{\mathrm{op}}}(T)$.
If this difference is empty, then $\mathfrak q\setminus\mathfrak p$ is nonempty
and the roles are reversed. This proves the vanishing.
\end{proof}

When $A=R$, the simple modules at $\mathfrak q$ and $\mathfrak p$
are $\kappa(\mathfrak q)$ and $\kappa(\mathfrak p)$, respectively.
Their injective envelopes are $I_R(R/\mathfrak q)$
and $I_R(R/\mathfrak p)$. Thus, Lemma~\ref{lem:closed-fiber} gives
\[
 \Tor_i^R
 \bigl(I_R(R/\mathfrak q),I_R(R/\mathfrak p)\bigr)=0
\]
for distinct $\mathfrak p,\mathfrak q\in\Spec R$ and every $i\ge0$.
This is the commutative vanishing result of Enochs and Jenda
\cite[Lemma~2.1]{EJ}. Enochs and Huang give the same central-action
argument for modules satisfying their property $t(\mathfrak p)$
\cite[Section~1(5)]{EH}. The proof above applies that argument
over a module-finite algebra.

\section{Auslander--Bridger formula}
Throughout this section, $A$ is an Iwanaga--Gorenstein 
Noether algebra. Its two
self-injective dimensions coincide by \cite[Lemma~A]{Zaks}; we denote
their common value by $\idim A$. 
We first use Tor nonvanishing to detect Gorenstein projective dimension
through injective envelopes of simple modules. Under the Auslander
condition, this detection result yields the Auslander--Bridger formula
in terms of the $S$-depths of $A$ and $M$.

We use the following standard facts; see \cite{AB,Holm,Iwanaga}.

\begin{lemma}\label{lem:ext-criterion}
Let $A$ be an Iwanaga--Gorenstein Noether algebra and let $M\in\mod A$ be nonzero.
The following statements hold.
\begin{enumerate}
    \item $\Gpdim_A M<\infty$.
    \item $\Gpdim_A M=\sup\{i\ge0\mid\Ext_A^i(M,A)\ne0\}$.
    \item $\Gpdim_A M=\pdim_A M$ whenever $\pdim_A M<\infty$.
    \item $\idim A=\sup\{\fdim_{A^{\op}}I\mid I \text{ an injective
    $A^{\op}$-module}\}$.
\end{enumerate}
\end{lemma}

We also use the fact that every Gorenstein projective module $G$ admits
a short exact sequence
\[
 0\lto G\lto P\lto G'\to0,
\]
where $P$ is projective and $G'$ is Gorenstein projective; see \cite{AB,Holm}.

The following bounds relate the grade of a simple module to the flat
dimension of its injective envelope; see \cite{KMT,Shen}.

\begin{lemma}\label{lem:weight}
Let $A$ be an Iwanaga--Gorenstein Noether algebra. For every simple
module $S$, we have
\[
 \grade_A S\le\fdim_A I_A(S).
\]
If $A$ is Auslander--Gorenstein, equality holds for all simple modules. 
The converse holds when $A$ is an Artin algebra.
\end{lemma}

The following lemma is the key step in the detection argument.

\begin{lemma}\label{lem:injective-detect}
Let $A$ be an Iwanaga--Gorenstein Noether algebra and 
let $M$ be nonzero in $\mod A$.
If
\[\Ext_A^p(M,A)\ne0,\]
then there exist a simple
$A^{\mathrm{op}}$-module $T$ and an integer $q\ge0$ such that
\[
 \Tor_{p+q}^A(I_{A^{\mathrm{op}}}(T),I^q(M))\ne0.
\]
In particular, $p+q\le\fdim_A I^q(M)$.
\end{lemma}
\begin{proof}
The $A^{\mathrm{op}}$-module $E=\Ext_A^p(M,A)$ is nonzero and finitely
generated, so it admits a simple quotient $T$.
Let $I=I_{A^{\mathrm{op}}}(T)$ and let $P_\bullet$ be a degreewise
finitely generated projective resolution of $M$. 
Finite projectivity gives an isomorphism 
\[
 I\otimes_A P_\bullet
 \cong \Hom_{A^{\mathrm{op}}}(\Hom_A(P_\bullet,A),I).
\]
Since $I$ is injective,
applying $\Hom_{A^{\mathrm{op}}}(-,I)$ is exact. Taking homology gives
\[
 \Tor_p^A(I,M)
 \cong \Hom_{A^{\mathrm{op}}}(\Ext_A^p(M,A),I).
\]
The composite $E\twoheadrightarrow T\rightarrowtail I$ is nonzero, 
so we have 
\[
 \Tor_p^A(I,M)\ne0.
\]

Set $Z^0=M$. The minimal injective coresolution of $M$ gives exact sequences
\[
 0\lto Z^q\lto I^q(M)\lto Z^{q+1}\lto0.
\]
If $\Tor_{p+q}^A(I,I^q(M))$ is zero, the connecting map
\[
 \Tor_{p+q+1}^A(I,Z^{q+1})
 \lto \Tor_{p+q}^A(I,Z^q)
\]
is surjective. If the indicated Tor group vanished for every
$q\geq 0$, these surjections would propagate the
nonvanishing to degree $\idim A+1$.
This contradicts the bound in Lemma~\ref{lem:ext-criterion}(4).
\end{proof}

We now show that one of the injective summands detected above is the
injective envelope of a simple module.

\begin{proposition}\label{prop:detect}
Let $A$ be an Iwanaga--Gorenstein Noether algebra and 
let $M$ be nonzero in $\mod A$. 
Then there exists a simple $A$-module $S$ such that
\[
 \Gpdim_A M\le\fdim_A I_A(S)-\depthr_S M.
\]
In particular, $\depthr_S M<\infty$.
\end{proposition}
\begin{proof}
Put $p=\Gpdim_A M$.
By Lemma~\ref{lem:ext-criterion}, $\Ext_A^p(M,A)\ne0$.
Lemma~\ref{lem:injective-detect} supplies a simple $A^{\op}$-module $T$
and an integer $q\ge0$ such that
\[
 \Tor_{p+q}^A(I_{A^{\op}}(T),I^q(M))\ne0.
\]
Write $I=I_{A^{\op}}(T)$.
An injective module over a Noether $R$-algebra decomposes as a direct sum of
indecomposable injectives. Since Tor commutes with direct sums,
some indecomposable summand $J$ of $I^q(M)$ satisfies
\[
 \Tor_{p+q}^A(I,J)\ne0.
\]
Choose $\mathfrak m\in\Max R$ in the support of the nonzero
finitely generated $R$-module $T$.
Nakayama's Lemma and simplicity give $\mathfrak mT=0$.
Then $T$ is also a simple module over $A_{\mathfrak m}^{\mathrm{op}}$.
By \cite[Remark~2.21]{KN}, injective envelopes agree under restriction,
so $I$ is isomorphic to $I_{A_{\mathfrak m}^{\mathrm{op}}}(T)$.

Write $J\cong I_{A_{\mathfrak p}}(S)$ by Lemma~\ref{lem:injective-classification}.
The nonzero Tor and Lemma~\ref{lem:closed-fiber} imply $\mathfrak p=\mathfrak m$.
The simple $A_{\mathfrak m}$-module $S$ is finite over $R_{\mathfrak m}$,
so Nakayama's Lemma and simplicity give $\mathfrak mS=0$.
The canonical isomorphism
\[
 A/\mathfrak mA\cong A_{\mathfrak m}/\mathfrak mA_{\mathfrak m}
\]
shows that $S$ is simple as an $A$-module via the localization map. 
Under this map, $J$ is isomorphic to $I_A(S)$. By
Lemma~\ref{lem:occurrence},
\[
 \depthr_S M+p\le q+p\le\fdim_A I_A(S).
\]
The claimed inequality and the finiteness of $\depthr_S M$ follow.
\end{proof}

We now assume that $A$ is Auslander--Gorenstein. Recall that this
means $A$ is Iwanaga--Gorenstein and $\fdim_A I^q(A)\le q$ for
every $q\ge0$ \cite{FGR}. An Auslander--Gorenstein algebra of finite
global dimension is called Auslander regular.

The next bound follows from \cite[Lemmas~3.1, 3.2 and 3.4]{HuangWeak}.
It is obtained by applying the Auslander condition for Gorenstein
projective modules to a finite Gorenstein projective resolution.

\begin{lemma}\label{lem:upper}
Let $A$ be an Auslander--Gorenstein Noether algebra and let $M$ be 
a nonzero $A$-module. Then
\[
 \fdim_A I^q(M)\le\Gpdim_A M+q
\]
for every $q\ge0$.
\end{lemma}

The following is the main result of this section.

\begin{theorem}\label{thm:ab}
Let $A$ be an Auslander--Gorenstein Noether algebra. For every nonzero
$M\in\mod A$, we have
\[
 \Gpdim_A M=\sup_{S\in\Sim A}
 \{\depthr_S A-\depthr_S M\}.
\]
In particular, $M$ is Gorenstein projective if and only if
\[
 \depthr_S M\ge\depthr_S A
\]
for every simple module $S$.
\end{theorem}
\begin{proof}
By Lemma~\ref{lem:weight}, 
we have $\grade_A S=\fdim_A I_A(S)$ for every simple module $S$.
It follows from Proposition~\ref{prop:detect} that 
\[
 \Gpdim_A M\le\sup_{S\in\Sim A}
 \{\depthr_S A-\depthr_S M\}.
\]

For the reverse bound, fix a simple module $S$ and 
put $g=\depthr_S A$. 
For every Gorenstein projective module $G$, 
take an exact sequence
\[
 0\lto G\lto P\lto G'\to0
\]
with $P$ projective and $G'$ Gorenstein projective. 
Since \(S\) admits a degreewise finitely generated projective resolution, 
the vanishing of $\Ext_A^i(S,A)$ for $0\le i<g$ 
extends to every projective module.
The embedding $G\rightarrowtail P$ gives
$\Hom_A(S,G)=0$ when $g>0$. 
For $1\le i<g$, the long exact
sequence gives
\[
 \Ext_A^i(S,G)\cong\Ext_A^{i-1}(S,G').
\]
Induction on $i$ gives the vanishing of $\Ext_A^i(S,G)$ for $0\le i<g$ 
and every Gorenstein projective module $G$.
Along a finite Gorenstein projective resolution of $M$ of length $p=\Gpdim_A M$,
dimension shifting gives
\[
 \Ext_A^i(S,M)=0\text{ for }0\le i<g-p.
\]
Then each term in the supremum is at most $p$, 
and equality holds for the simple module supplied by 
Proposition~\ref{prop:detect}.
\end{proof}

For a commutative Gorenstein ring $R$ of finite Krull dimension,
the last assertion of Theorem~\ref{thm:ab} recovers the characterization
of finitely generated Gorenstein projective modules as maximal
Cohen--Macaulay modules.

Indeed, the inequalities in Theorem~\ref{thm:ab} are equivalent to
\[
 \depthr_{R_{\mathfrak m}} M_{\mathfrak m}
 =\depthr R_{\mathfrak m}
\]
for every maximal ideal $\mathfrak m$ in the support of $M$.
Since each $R_{\mathfrak m}$ is Gorenstein and hence Cohen--Macaulay,
its depth equals its Krull dimension.

Over an Auslander regular algebra, every module has finite projective
dimension, so Theorem~\ref{thm:ab} gives the following Auslander--Buchsbaum formula.

\begin{corollary}\label{cor:ab}
Let $A$ be an Auslander regular Noether algebra. For every nonzero
$M\in\mod A$, we have
\[
 \pdim_A M=\sup_{S\in\Sim A}
 \{\depthr_S A-\depthr_S M\}.
\]
\end{corollary}

\begin{remark}\label{rem:weaker-hypothesis}
Theorem~\ref{thm:ab} remains valid if $A$ is only assumed to be an 
Iwanaga--Gorenstein Noether algebra satisfying
\[
 \grade_A S=\fdim_A I_A(S)
\]
for every simple $A$-module $S$. Indeed, the proof uses the
Auslander condition only through this equality. The detection argument
and the dimension shifting apply under the stated hypothesis.
\end{remark}

\begin{example}\label{ex:unequal-grades}
Let $k$ be a field and let $Q$ be the quiver
\[
 \begin{tikzpicture}[baseline=(one.base),>=stealth]
  \node (one) at (0,0) {$1$};
  \node (two) at (1.4,0) {$2$};
  \draw[->] (one) to[bend left=24] node[midway,above] {$\alpha$} (two);
  \draw[->] (two) to[bend left=24] node[midway,below] {$\beta$} (one);
 \end{tikzpicture}
 .
\]
Let $A=kQ/\langle\beta\alpha\rangle$, where
$\beta\alpha$ means $\alpha$ followed by $\beta$.
Write $S_i$ for the vertex simple modules
and $P_i$ for their projective covers.

The composition factors of
$P_1$ and $P_2$, from top to socle, are respectively
$(S_1,S_2)$ and $(S_2,S_1,S_2)$. Then
$\operatorname{rad}P_1\cong S_2$ and
$\operatorname{rad}P_2\cong P_1$, giving
$\pdim_A S_2=1$ and $\pdim_A S_1=2$. In particular, $A$ has global
dimension two.

The paths ending at each vertex give isomorphisms
\[I_A(S_2)\cong P_2 \text{ and } I_A(S_1)\cong P_2/\soc P_2.\] 
Since $P_2/P_1\cong S_2$,
the minimal injective coresolution of $P_1$
is
\[
 0\lto P_1\lto P_2\lto P_2
 \lto I_A(S_1)\lto0.
\]
As $P_2$ is injective, the first two terms of the minimal injective
coresolution of $A$ are projective. The global dimension
bound gives the Auslander condition in degree two. By
Lemma~\ref{lem:occurrence}, the first occurrences of $I_A(S_1)$ and
$I_A(S_2)$ give $\grade_A S_1=2$ and $\grade_A S_2=0$.

Take $M=S_2$. The nonsplit exact sequence 
\[0\lto S_2\lto P_1\lto S_1\lto0\]
and the vanishing of
$\Hom_A(S_1,S_2)$ give $\depthr_{S_1} M=1$, while
$\depthr_{S_2} M$ is zero.
Substituting these values into Theorem~\ref{thm:ab} gives
\[
 \Gpdim_A M=\sup\{2-1,\,0-0\}=1.
\]
The supremum is attained by $S_1$.
This agrees with the projective resolution
\[0\lto P_1\lto P_2\lto M\lto0.\]
\end{example}

Although a simple module outside the socle can attain the supremum,
maximal Gorenstein projective dimension is detected by a simple submodule.

\begin{corollary}\label{cor:socle-detection}
Let $A$ be an Auslander--Gorenstein Noether algebra and
let $M$ be nonzero in $\mod A$. The following are equivalent.
\begin{enumerate}
\item $\Gpdim_A M=\idim A$;
\item $\fdim_A I_A(M)=\idim A$;
\item $M$ contains a simple submodule $S$ with $\grade_A S=\idim A$.
\end{enumerate}
\end{corollary}
\begin{proof}
For every simple module $S$, Lemma~\ref{lem:weight} gives
\[
 \grade_A S=\fdim_A I_A(S)\le\idim A.
\]
A term in the supremum in Theorem~\ref{thm:ab} equals $\idim A$
if and only if
\[
 \grade_A S=\idim A \text{ and }\depthr_S M=0.
\]
Since the supremum in Theorem~\ref{thm:ab} is attained,
(1) and (3) are equivalent.

If (3) holds, $I_A(S)$ is a direct summand of $I_A(M)$.
Every injective module has flat dimension at most $\idim A$,
so (2) follows. Conversely, Lemma~\ref{lem:upper} gives
\[
 \fdim_A I_A(M)\le\Gpdim_A M\le\idim A,
\]
which proves that (2) implies (1).
\end{proof}

The Auslander condition in Corollary~\ref{cor:socle-detection} 
cannot simply be omitted.

\begin{example}
Let $k$ be a field and let $A$ be the path algebra of the quiver
\[
 \begin{tikzpicture}[baseline=(two.base),>=stealth]
  \node (one) at (-1.2,0) {$1$};
  \node (two) at (0,0) {$2$};
  \node (three) at (1.2,0) {$3$};
  \draw[->] (two) -- node[midway,above] {$\alpha$} (one);
  \draw[->] (two) -- node[midway,above] {$\beta$} (three);
\end{tikzpicture}
.
\]
This algebra is hereditary. Write $S_i$ for its vertex
simple modules, $P_i$ for their projective covers, and put $M=I_A(S_1)$.
There is an exact sequence
\[
 0\lto P_3\lto P_2\lto M\lto0.
\]
Then $\Gpdim_A M=\pdim_A M=1$. On the other hand, $\soc M\cong S_1$, and $S_1$ is
projective, so $\grade_A S_1=0$. Every simple submodule $S$ of $M$ therefore satisfies
$\grade_A S=0$, so condition (3) is not satisfied. Since $M$ is injective with socle
$S_1$, only $S_1$
contributes a finite term to the supremum in the proposed
formula. That supremum is zero, whereas $\Gpdim_A M=1$.
\end{example}

We finish with a converse in the Artin setting.
Recall that an Artin algebra is an algebra over a commutative Artinian ring
that is finitely generated as a module over that ring.
Let $A$ be a Noether $R$-algebra.
Then $A$ has finite length as an $R$-module if and only if it
is an Artin algebra over $R/\ann_RA$.

\begin{corollary}\label{cor:artin-characterization}
Let $A$ be an Iwanaga--Gorenstein Artin algebra.
The following are equivalent.
\begin{enumerate}
\item $A$ is Auslander--Gorenstein;
\item The Auslander--Bridger formula holds for every nonzero
$M\in\mod A$.
\end{enumerate}
\end{corollary}
\begin{proof}
Theorem~\ref{thm:ab} shows that (1) implies (2).
Conversely, assume (2) and let $S$ be a simple $A$-module.
Taking $M=I_A(S)$, injectivity and the equality
$\soc I_A(S)=S$ show that only $S$ contributes a finite term
to the supremum. Then 
\[
 \Gpdim_A I_A(S)=\grade_A S.
\]
By Lemma~\ref{lem:ext-criterion}(4), the module $I_A(S)$ 
has finite flat dimension. Since flat modules over an 
Artin algebra are projective, Lemma~\ref{lem:ext-criterion} gives
\[
 \pdim_A I_A(S)=\Gpdim_A I_A(S)=\grade_A S.
\]
Lemma~\ref{lem:weight} now shows that $A$ is Auslander--Gorenstein.
\end{proof}

The simple-module characterization used in the proof appears in
\cite{KMT,Shen}. The forward implication can also be obtained
from Huang's characterization of Gorenstein projective dimension
by flat-dimension bounds on minimal injective coresolutions
\cite{HuangWeak}.
Over an Artin algebra, each injective term is a direct sum of
injective envelopes of simple modules. Reading these bounds at
the first occurrence of each $I_A(S)$, and applying
Lemmas~\ref{lem:occurrence} and~\ref{lem:weight}, recovers
the Auslander--Bridger formula.

\section{The defect criterion}
For an Artin algebra $A$ and $M\in\mod A$,
the dominant dimension of $M$ is the first degree of the minimal
injective coresolution of $M$ whose term is nonprojective,
with value $\infty$ if every term is projective.
An algebra $A$ is called minimal Auslander--Gorenstein if its
self-injective dimension is finite and at most its dominant dimension.
A minimal Auslander--Gorenstein algebra of finite global dimension
is called a higher Auslander algebra.
The defect $\operatorname{def} A$ is the number of isomorphism classes
of nonprojective indecomposable injective $A$-modules.

Let $A$ be an Iwanaga--Gorenstein Artin algebra and put $d=\idim A$.
For $d=0$, $A$ is self-injective.
Then $\domdim_A A=\infty$ and $\operatorname{def} A$ is zero,
while every simple module has Gorenstein projective dimension zero.
For $d\ge1$, the socle criterion of
Li, Marczinzik and Zhang \cite{LMZ} says that $A$ is minimal
Auslander--Gorenstein precisely when, for every simple $S$,
\[
 I_A(S)\text{ is projective}\Longleftrightarrow\,\Gpdim_A S<d.
\]

Counting the complementary simple modules gives the following theorem.

\begin{theorem}\label{thm:defect}
Let $A$ be an Iwanaga--Gorenstein Artin algebra of positive self-injective
dimension $d$. Then
\[
 \#\{S\in\Sim A\mid\Gpdim_A S=d\}
 \ \le\ \operatorname{def} A.
\]
Equality holds if and only if $A$ is minimal Auslander--Gorenstein.
\end{theorem}
\begin{proof}
If $I_A(S)$ is projective, the envelope sequence gives
\[
 \Ext_A^d(S,A)\cong\Ext_A^{d+1}(I_A(S)/S,A)=0.
\]
Lemma~\ref{lem:ext-criterion} gives $\Gpdim_A S<d$.
The assignment $S\mapsto I_A(S)$ defines an injection from the isomorphism classes
of simple modules of Gorenstein projective dimension $d$ to those of
nonprojective indecomposable injective modules. Both sets are
finite, so equality holds exactly when the injection is surjective.
By the socle criterion of Li, Marczinzik and Zhang, this is equivalent
to $A$ being minimal Auslander--Gorenstein.
\end{proof}

We obtain the following corollary; see also
\cite[Corollary~5.7(3)]{MarFin}.

\begin{corollary}\label{cor:defect-one}
Every Iwanaga--Gorenstein Artin algebra of defect one
is minimal Auslander--Gorenstein.
\end{corollary}
\begin{proof}
Since a self-injective Artin algebra has defect zero, $d=\idim A>0$.
Choose a simple submodule $S$ of the nonzero last term $I^d(A)$.
By Lemma~\ref{lem:occurrence}, $\Ext_A^d(S,A)$ is nonzero and 
$\Gpdim_A S=d$ by Lemma~\ref{lem:ext-criterion}. The count in
Theorem~\ref{thm:defect} is therefore at least one, so equality holds.
\end{proof}

\begin{corollary}\label{cor:higher-defect}
Let $A$ be an Artin algebra of finite global dimension $d\ge1$.
Then $A$ is higher Auslander if and only if
\[
 \operatorname{def} A=\#\{S\in\Sim A\mid\pdim_A S=d\}.
\]
In particular, if $\operatorname{def} A=1$, then $A$ is higher Auslander.
\end{corollary}
\begin{proof}
Since $A$ has finite global dimension, 
$\Gpdim_A S=\pdim_A S$ for every simple module $S$.
The equivalence follows from Theorem~\ref{thm:defect}, and the
last assertion follows from Corollary~\ref{cor:defect-one}.
\end{proof}

For $n\geq 2$, the connected linear Nakayama algebra 
of defect one with $n$ simple modules has Kupisch series 
$((2)^{n-1},1)$. Its global and dominant dimensions 
both equal $n-1$; see \cite[Section~3]{MMZ}.

\section{Comparison with known formulas}\label{sec:comparison}

In this section, we derive a dominant dimension formula for modules over
minimal Auslander--Gorenstein Artin algebras. We then express the main
theorem as a depth identity and give its specialization to commutative
Gorenstein rings.

For an Artin algebra $A$, put
\[
 \mathcal{U}=\{S\in\Sim A\mid I_A(S)\text{ is not projective}\}.
\]
For a finitely generated module $M$, its dominant dimension is given by
\[\domdim_AM=\inf_{S\in \mathcal U}\{\depthr_S M\}.\]

\begin{corollary}\label{cor:min}
Let $A$ be a minimal Auslander--Gorenstein Artin algebra and let
$M\in\mod A$. If $M\oplus\Sigma M$ is not Gorenstein projective, then
\[
 \Gpdim_A M+\domdim_A M=\idim A.
\]
Here $\Sigma M$ denotes the minimal cosyzygy of $M$.
\end{corollary}
\begin{proof}
Put $d=\idim A$. For $q<d$, the term $I^q(A)$ is projective, 
and $I^d(A)$ has projective dimension at most $d$. 
Then $A$ is Auslander--Gorenstein. By Lemmas~\ref{lem:occurrence}
and~\ref{lem:weight}, the modules in $\mathcal{U}$ have grade $\idim A$,
while the others have grade zero. If $\Gpdim_A M>0$, the supremum in
Theorem~\ref{thm:ab} is attained in $\mathcal{U}$. Taking the infimum 
of the $S$-depths of $M$ over $S\in\mathcal U$ gives $\domdim_A M$.

If $M$ is Gorenstein projective, then $\Sigma M$ is not. 
By Theorem~\ref{thm:ab} we have 
\[\depthr_S M\ge\idim A \text{ for } S\in \mathcal{U}.\] 
It follows that 
\[\domdim_A M\ge\idim A\geq 1.\] 
Then $I^0(M)$ is projective and the exact sequence
\[
 0\lto M\lto I^0(M)\lto\Sigma M\lto0
\]
gives $\Gpdim_A\Sigma M=1$. Applying the first case to $\Sigma M$ yields
\[
 \domdim_A M=1+\domdim_A\Sigma M=\idim A.
\]
This completes the proof.
\end{proof}

Under the hypotheses of \cite[Theorem~2.2]{MarMin}, 
Corollary~\ref{cor:min} also follow from the characterization
of Gorenstein projective dimension by dominant dimension.
The proof above also treats Gorenstein projective modules whose
cosyzygies are not Gorenstein projective. Related dimension formulas
are discussed in \cite{ChIM}. For a higher Auslander algebra of global dimension $d$,
Gorenstein projective dimension agrees with projective dimension. If $M$ is projective
and noninjective,
its envelope sequence is nonsplit, so $\Sigma M$ is nonprojective.
The corollary yields an equality 
\[\pdim_A M+\domdim_A M=d\]
for every module $M$ that is not projective-injective; see \cite{CM}.

Let $A$ be an Iwanaga--Gorenstein Noether algebra. 
For nonzero $M$ in $\mod A$, define
\[
 \depthr_A M=\inf_{S\in\Sim A}
 \{\idim A-\depthr_S A+\depthr_S M\}.
\]
If $A$ is Auslander--Gorenstein, 
Theorem~\ref{thm:ab} rewrites as
\[
 \Gpdim_A M+\depthr_A M=\idim A.
\]
Moreover, if all simple modules have the same grade, then 
\[
 \depthr_A M=\inf_{S\in\Sim A}
 \{\depthr_S M\}.
\]
When the hypotheses of both results hold, the resulting formula
for modules of finite projective dimension agrees with that of
Wu and Zhang \cite[Theorem~0.1]{WZ}.

Let $R$ be a commutative Gorenstein ring of finite Krull dimension.
Then $R$ is Auslander--Gorenstein \cite{Bass}.
For each maximal ideal $\mathfrak m$, $R_{\mathfrak m}$ is Gorenstein, so 
\[
 \depthr R_{\mathfrak m}=\idim R_{\mathfrak m}
 =\dim R_{\mathfrak m}=\operatorname{ht}\mathfrak m.
\]
For every $N$ in $\mod R$, localization gives
\[
 \inf\{i\ge 0\mid \Ext_R^i(R/\mathfrak m,N)\ne 0\}
 =\depthr_{R_{\mathfrak m}} N_{\mathfrak m},
\]
Applying this identity to $R$ and $M$, Theorem~\ref{thm:ab}
gives
\[
 \Gpdim_R M
 =\sup_{\mathfrak m\in\Max R}
 \bigl\{\depthr R_{\mathfrak m}
       -\depthr_{R_{\mathfrak m}} M_{\mathfrak m}\bigr\}
\]
for nonzero $M$ in $\mod R$.
When $R$ is local, this is the Auslander--Bridger formula
\[
 \Gpdim_R M+\depthr_R M=\depthr R;
\]
see \cite[Theorem~4.13]{AB}.
For $M$ of finite projective dimension, this gives the
classical Auslander--Buchsbaum formula;
see \cite[Theorem~3.7]{AusBuch}.

\section*{Acknowledgements}
This work is partially supported by the National Natural Science Foundation of
China (Nos.~12571036 and 11801141).

\section*{Declaration of AI use}
The author used ChatGPT (GPT-6) to assist with language editing,
exposition, and preliminary checks of mathematical arguments.
All suggestions were reviewed and verified by the author, who takes full
responsibility for the final manuscript.

\end{document}